\documentclass[a4paper,10pt,twoside]{scrartcl}
\usepackage{german,amssymb,amsmath,latexsym,epsfig,graphics,graphicx,mathrsfs,color,amsthm,cite}
\usepackage[latin1]{inputenc}
\usepackage{makeidx}

\newtheorem{theorem}{Theorem}[section]

\newtheorem{definition}[theorem]{Definition}
\newtheorem{lemma}[theorem]{Lemma}

\newtheorem*{pf}{Proof}

\numberwithin{equation}{section}

\newcommand{\R}{{\mathbb R}}
\newcommand{\N}{{\mathbb N}}

\DeclareMathOperator*{\llim}{\lambda-lim}

\begin{document}
\title{\Large{Notes on Spaces of Functions Converging at Infinity}}

\author{\large{Nico Tauchnitz} \\ {\small Nico.Tauchnitz@hszg.de}}
\date{\large{\today}}
\maketitle

\begin{abstract}
\textbf{Abstract.} This preprint concerns Banach spaces of functions converging at infinity.
In particular, spaces of continuous functions, Lebesgue spaces and sequence spaces.
In each framework we show versions of Riesz's representation theorem. \\[1mm]
\end{abstract}


\section{Motivation} \label{SectionMotivation}
The need of continuous functions converging at infinity arose in our considerations of infinite horizon optimal control
problems in \cite{Tauchnitz}.
In this framework, the methods of convex optimization led to challenges hidden deep inside the complex theory. 
As an illustration, consider the constraint
$$K(t)=\int_0^t k(s) \,ds, \qquad k(t) \geq 0, \quad K(t) \leq K \mbox{ for all } t \geq 0.$$
Then the continuous trajectory $x(\cdot)=K(\cdot)-K$ satisfies the inequality $x(t)\leq 0$ on $[0,\infty)$ iff
\begin{enumerate}
\item $x(\cdot)$ belongs to the cone
      $\mathscr{K}=\{ z(\cdot) \in C([0,\infty),\R) \,|\, z(t)\leq 0 \mbox{ for all } t \geq 0\}$;
\item the inequality $f\big(x(\cdot)\big)=\sup\limits_{t \geq 0} x(t)\leq 0$ holds.
      Thereby, $f$ is convex on $C([0,\infty),\R)$.
\end{enumerate}
The particular space $C_0([0,\infty),\R)$ of continuous functions vanishing at infinity
provides a Riesz' representation theorem including the regular signed Borel measures on $[0,\infty)$ (cf. \cite{Rudin}).
But in this framework the following degenerations arising:
\begin{enumerate}
\item The cone $\mathscr{K}=\{ z(\cdot) \in C_0([0,\infty),\R) \,|\, z(t)\leq 0 \mbox{ for all } t \geq 0\}$
      possesses an empty interior:
      Let $z(\cdot) \in \mathscr{K}$ with $z(t) <0$ on $[0,\infty)$ and consider the sequence
      $$z_n(t)=\left\{\begin{array}{ll} z(t), & t \not\in [2n-1,2n+1], \\
                                        z(t)+2|z(2n)| \cdot (t-(2n-1))((2n+1)-t), & t \in [2n-1,2n+1].
                      \end{array} \right.$$
      Then $\|z_n(\cdot)-z(\cdot)\|_\infty \to 0$ as $n \to \infty$ and $z_n(\cdot) \not\in \mathscr{K}$ for all $n$.
      Hence in the space $C_0([0,\infty),\R)$ the interior of $\mathscr{K}$ and the origin cannot be separated by the
      Hahn-Banach separation theorem for convex sets.
      Therefore, the basic tool in convex optimization cannot be applied. 
\item In the space $C_0([0,\infty),\R)$ the subdifferential of
      $f\big(x(\cdot)\big)=\sup\limits_{t \geq 0} x(t)$
      delivers in $x(\cdot)=0$ the set of all non-negative regular Borel measures $\mu$ on $[0,\infty)$ with the
      non-strict attribute $\|\mu\| \leq 1$. \\
      Indeed, by definition of $\partial f(0)$ we obtain in the space $C_0(\R_+,\R)$ the inequalities
      $$\sup_{t \in \R_+} x(t) \geq \int_0^\infty x(t) \, d\mu(t)
        \geq - \sup_{t \in \R_+} \big( -x(t)\big)
        = \inf_{t \in \R_+} x(t).$$
      This shows $\mu \geq 0$ and $\|\mu\| \leq 1$.
      Moreover, for any function $z(\cdot) \in C_0([0,\infty),\R)$ with $z(t) <0$ on $[0,\infty)$ the subdifferential   
      $\partial f\big(z(\cdot)\big) =\{x^* \in \partial f(0) \,|\, f\big(z(\cdot)\big)= \langle x^*,z(\cdot) \rangle \}$
      delivers
      $$0 = \sup_{t \geq 0} z(t) =\int_0^\infty z(t) \, d\mu(t) \quad\Leftrightarrow\quad \mu=0.$$
      But $\|\mu\|=0$ may lead to the degenerated case of trivial Lagrange multipliers. 
\end{enumerate}
In \cite{Tauchnitz} we proposed to consider the space of continuous functions converging at infinity.
Then the optimization methods require a Riesz' representation theorem.
As a path to Riesz' representation theorem in the space of continuous functions converging at infinity one may use a
compactification of $[0,\infty]$, e.g. by $t(s)=s/(1+s)$, $s \in [0,\infty]$.
The inverse mapping $s(t)= t/(1-t)$ maps $[0,1]$ onto $[0,\infty]$.
Then the transformation generating the additional weight $1/(1+s)^2$. 
Therefore, a Riesz's representation theorem achieved by transformation of $[0,1]$ onto $[0,\infty]$ depends on the
particularly chosen transformation.
In the subsequent sections we provide a transformation free statement of the Riesz's representation theorem.
Our approach follows the methods of functional analysis.


\section{Basic Notations and Definitions}
Throughout this preprint we make frequently use of the notations
$$\R_+=[0,\infty), \quad \overline{\R}_+=[0,\infty], \quad \R=(-\infty,\infty), \quad \overline{\R}=[-\infty,\infty].$$
Thereby,
$\overline{\R}=[-\infty,\infty]$ denotes the
affinely extended real number system obtained from the real number system $\R$ by adding the elements
$\infty$ and $-\infty$.
By the continuous mapping $\displaystyle x \to \frac{x}{1+|x|}$
the extended real number line $\overline{\R}$ is homeomorphic to the interval $[-1,1]$.
Hence, $\overline{\R}$ is an open and compact subset of $\overline{\R}$.
A set $U$ is a neighborhood of $\infty$, if and only if it contains a set $\{x \,|\, x > a\}$ for some real number $a$. 
Neighborhoods of $-\infty$ can be defined analogously.

\begin{definition} \label{DefinitionBorel}
\begin{enumerate}
\item[(a)] Let $I \subseteq \R$.
           The Borel algebra $\mathscr{B}(I)$ on $I$ consists of all open subsets of $I$.
\item[(b)] The Borel algebra $\mathscr{B}(\overline{\R}_+)$ on $\overline{\R}_+$ consists of all sets $B=A \cup E$ with
           $A \in \mathscr{B}(\R_+)$ and $E \subseteq \{\infty\}$.
\item[(c)] The Borel algebra $\mathscr{B}(\overline{\R})$ on $\overline{\R}$ consists of all sets $B=A \cup E$ with
           $A \in \mathscr{B}(\R)$ and $E \subseteq \{-\infty,\infty\}$.
\end{enumerate}
\end{definition}

\begin{definition}
A Borel measure $\mu$ defined on the Borel algebra $\mathscr{B}(I)$, $I \subseteq \overline{\R}$, is regular if
$$\mu(A) = \sup\{ \mu (B) \,|\, B\subseteq A, \; B \mbox{ compact and measurable}\}$$
and
$$\mu (A)=\inf\{\mu(C) \,|\, C\supseteq A, \; C \mbox{open and measurable}\}$$
for any $A \in \mathscr{B}(I)$.
\end{definition}

Let $\mu$ be a signed measure on $\overline{\R}$ and let the system $\{E_i\}$ of measureable sets
$E_i$ be a countable partition of $\overline{\R}$.
As stated in Rudin \cite{Rudin} in the case of $\R_+$, the absolute convergence of the series
$\displaystyle \sum_{i=1}^\infty \mu(E_i)$ is now part of the requirements of a signed measure on $\overline{\R}$.
Let $\nu$ be a signed regular Borel measure on $\mathscr{B}(\R)$,
let $\nu_{-\infty}$ and $\nu_\infty$ are signed measures concentrated at $t=-\infty$ and $t=\infty$, respectively.
Since $\overline{\R}$ is open and compact,
a signed Borel measure $\mu$ on $\mathscr{B}(\overline{\R})$ is regular
if and only if $\mu$ possesses the representation $\mu=\nu+\nu_{-\infty}+\nu_\infty$. \\

\begin{definition} \label{DefinitionBorelmass}
\begin{enumerate}
\item[(a)] Let $I \subseteq \R$.
           $\mathscr{M}(I)$ denotes the set of signed regular Borel measures on $\mathscr{B}(I)$.
\item[(b)] By $\mathscr{M}(\overline{\R}_+)$ we denote the set of the signed regular Borel measures on
           $\mathscr{B}(\overline{\R}_+)$.
           Any $\mu \in \mathscr{M}(\overline{\R}_+)$ has the unique representation $\mu=\nu+\nu_\infty$,
           where $\nu \in \mathscr{M}(\R_+)$ and $\nu_\infty$ is a signed regular Borel measure concentrated at $t =\infty$.
\item[(c)] By $\mathscr{M}(\overline{\R})$ we denote the set of the signed regular Borel measures on
           $\mathscr{B}(\overline{\R})$.
           Any $\mu \in \mathscr{M}(\overline{\R})$ has the unique representation $\mu=\nu+\nu_{-\infty}+\nu_\infty$,
           where $\nu \in \mathscr{M}(\R)$ and where
           $\nu_{-\infty}$ and $\nu_\infty$ are signed regular Borel measure concentrated at $t=-\infty$ and $t=\infty$,
           respectively.
\end{enumerate}
\end{definition}

In the Definitions \ref{DefinitionBorel} and \ref{DefinitionBorelmass} we introduced the objects:
\begin{center}\begin{tabular}{lcl}
$\mathscr{B}(I)$ &--& The Borel algebra on $I \subseteq \overline{\R}$. \\[1mm]
$\mathscr{M}(I)$ &--& The set of regular signed Borel measures on $\mathscr{B}(I)$ \hspace*{2.3cm}. 
\end{tabular}\end{center}

In the subsequent sections the following Banach spaces are of our interest:
\begin{center}\begin{tabular}{lcl}
$C_{\lim}(\R_+,\R^n)$ &--& The space of continuous functions converging at infinity (Section \ref{SectionContinuous}). \\[1mm]
$C_{\lim}(\R,\R^n)$ &--& The space of continuous functions converging at
                         $t=\pm \infty$ (Section \ref{SectionContinuous}). \\[1mm]
$L_{p,\lim}(\R_+,\R^n)$ &--& Lebesgue spaces of functions converging in measure at infinity 
                             (Section \ref{SectionLebesgue}). \\[1mm]
$L_{p,\lim}(\R,\R^n)$ &--& Lebesgue spaces of functions converging in measure at $t=\pm \infty$
                           (Section \ref{SectionLebesgue}). \\[1mm]
$W^1_{p,\lim}(\R_+,\R^n)$ &--& Sobolev spaces of functions converging in measure at infinity
                              (Section \ref{SectionRemarks}). \\[1mm]
$\ell_{p,\lim}$ &--& Spaces of convergent sequences (Section \ref{SectionRemarks}).
\end{tabular}\end{center}


\begin{definition} \label{DefinitionConvergence}
\begin{enumerate}
\item[(a)] We say $x(\cdot):\R_+ \to \R^n$ converges at infinity to $a \in \R^n$
           if for every $\varepsilon >0$ there exists $N>0$ such that
           $\|x(t)-a\|\leq \varepsilon$ for all $t \geq N$.
           Then we write $\lim\limits_{t \to \infty} x(t)=a$.
\item[(b)] We say $x(\cdot):\R_+ \to \R^n$ converges in Lebesgue measure $\lambda$ at infinity to $a \in \R^n$
           if for every $\varepsilon >0$,
           $\lim\limits_{N \to \infty} \lambda\big(\{t\geq N \,|\, \|x(t)-a\|\geq \varepsilon\}\big)=0$.
           The $\lambda-$limit is denoted by $\llim\limits_{t \to \infty} x(t)$.
\end{enumerate}
\end{definition}

%

\begin{definition} \label{DefinitionMappings}
Let $X,Y$ be Banach spaces, let $\Lambda:X \to Y$ be linear and continuous.
\begin{enumerate}
\item[(a)] The set $Im\,\Lambda = \{y \in Y \,|\, \Lambda x=y,\; x\in X\}$ denotes the range of $\Lambda$.
\item[(b)] The set $Ker\,\Lambda = \{x \in X \,|\, \Lambda x=0\}$ is the kernel of $\Lambda$.
\item[(c)] $(Ker\,\Lambda)^\perp=\{x^* \in X^*\,|\, \langle x^*,x \rangle=0 \mbox{ for all } x \in Ker\,\Lambda\}$
           denotes the annihilator of $Ker\,\Lambda$.
\item[(d)] The adjoint operator $\Lambda^*:Y^*\to X^*$ of $\Lambda$ is the linear operator defined by
           $\langle y^*,\Lambda x \rangle=\langle \Lambda^* y^*, x \rangle$
           for all $x \in X$ and $y^* \in Y^*$.
\item[(e)] The set $Im\,\Lambda^*=\{x^* \in X^*\,|\, x^*=\Lambda^* y^* , y^* \in Y^*\}$ denotes the range of $\Lambda^*$.
\end{enumerate}
\end{definition}
\section{Basic Concepts} \label{SectionConcepts}
In the subsequent sections we achieve different versions of Riesz' representation theorem in spaces of
functions converging at infinity.
Our basic tool is the Closed range theorem (cf. \cite{Werner}).

\begin{theorem}[Closed Range Theorem] \label{LemmaClosedRange}
Let $X,Y$ be Banach spaces, let $\Lambda:X \to Y$ be linear and continuous.
Then $\Lambda$ has a closed range in $Y$ if and only if $(Ker\,\Lambda)^\perp=Im\,\Lambda^*$.
\end{theorem}

In this section we outline the specific concepts, where
the first concept concerns the functions on $\R_+$, while the second concept aims continuous functions on $\R$.
Throughout this section let
\begin{itemize}
\item[$\cdot$] $X_0$ be a given Banach space of mappings $x_0(\cdot):\R_+ \to \R^n$ with norm $\|\cdot\|_{X_0}$.
\item[$\cdot$] $X$ be the space of all $x(\cdot):\R_+ \to \R^n$,
               which possesses a unique representation $x(t)=x_0(t)+a$ on $\R_+$, where $x_0(\cdot) \in X_0$ and $a \in \R^n$.
\item[$\cdot$] $X$ be equipped with a $p-$norm,
               e.g. $\|x(\cdot)\|^p_X=\|x_0(\cdot)\|^p_{X_0}+ \|a\|^p$ if $1\leq p <\infty$ and
               $\|x(\cdot)\|_X=\|x_0(\cdot)\|_{X_0}+ \|a\|$ if $p=\infty$.
               Then $X$ becomes a Banach space.
\end{itemize}

{\bfseries First scheme:}
The subsequent result provides directly the versions of Riesz' representation theorem in Banach spaces $X$,
where $X_0$ is $C_0(\R_+,\R^n)$, $L_p(\R_+,\R^n)$ or $\ell_p$ with $1\leq p<\infty$.

\begin{lemma} \label{LemmaRepresentation1}
Any continuous linear functional $x^*$ on $X$ can be uniquely represented in the form
$$\langle x^*(\cdot),x(\cdot) \rangle = \langle x_0^*(\cdot),x(\cdot)-a \rangle + \alpha^T a,
  \quad \|x^*\|=\|x_0^*\|+\|\alpha\|, \quad x_0^* \in X_0^*,\; \alpha \in \R^n.$$
\end{lemma}

\textbf{Proof.}
Any $x(\cdot) \in X$ has the unique representation $x(\cdot)=x_0(\cdot)+a$ with $x_0(\cdot) \in X_0$ and $a \in \R^n$.
Consider the linear mapping $\Lambda$ with $\Lambda x(\cdot) = x_0(\cdot)$, which maps $X$ onto $X_0$.
Then $Ker\,\Lambda$ consists of all $x(\cdot) \in X$ with $x_0(\cdot)=0$.
By definition of $\|\cdot\|_X$, $\Lambda$ is continuous.
Let $x^* \in X^*$.
We denote by $\alpha_i$ the value of the functional $x^*$ at the vector-valued function whose $i$-th component
is identical $1$ and whose remaining components are identically zero.
Now, we consider the functional $x_1^* \in X^*$ defined by the formula
$$\langle x_1^*,x(\cdot) \rangle = \langle x^*,x(\cdot) \rangle - \alpha^T a,
  \quad x(t)=x_0(t)+a, \quad \alpha=(\alpha_1,...,\alpha_n).$$
Obviously, $x_1^* \in (Ker\,\Lambda)^\perp$.
By Theorem \ref{LemmaClosedRange},
there exists a functional $x_0^* \in X^*_0$ with $x_1^* = \Lambda^* x_0^*$.
That means the equation
$\langle x_1^*,x(\cdot) \rangle = \langle \Lambda^* x_0^*,x(\cdot) \rangle 
                                 = \langle x_0^*,\Lambda x(\cdot)\rangle =\langle x_0^*,x_0(\cdot) \rangle$
holds for all $x(\cdot) \in X$.
Consequently, for $x^* \in X^*$ we obtain the representation
$$\langle x^*,x(\cdot) \rangle = \langle x_1^*,x(\cdot) \rangle + \alpha^T a
  = \langle x_0^*, x_0(\cdot) \rangle + \alpha^T a= \langle x_0^*, x(\cdot)-a \rangle + \alpha^T a.$$
Finally, the uniqueness of this representation can be verified directly. \hfill $\square$ \\

{\bfseries Second scheme:}
Our consideration of functions $x(\cdot)$ on $\R$ bases on a splitting argument.
Each $x(\cdot)$ can be uniquely split into the parts on the positive and negative real axis by
$$\big(x_1(t),x_2(t)\big)=\big(x(-t),x(t)\big), \qquad t \in \R_+.$$
In the particular case of continuous functions we have to preserve the continuity at $t=0$ and determine the splitting by
$$\big(x_1(t),x_2(t)\big)=\big(x(-t)-x(0),x(t)-x(0)\big), \qquad t \in \R_+.$$
According to the latter splitting, we introduce the subspace $Y=\{y(\cdot)\in X \,|\, y(0)=0\}$ of $X$.
Therefore, any $y(\cdot) \in Y$ possesses the representation $y(t)=x_0(t)+a$ with $x_0(\cdot) \in X_0$, $a \in \R^n$
and with the additional attribute $y(0)=0$.
Then the splitting leads to a mapping, which maps onto the product $Y \times Y$.
Thus, we have to provide a representation formula on $Y$:

\begin{lemma} \label{LemmaRepresentation2}
Any continuous linear functional $y^*$ on $Y$ can be uniquely represented in the form
$$\langle y^*(\cdot),y(\cdot) \rangle = \langle x_0^*,y(\cdot)-a \rangle,
  \quad \|y^*\|=\|x_0^*\|, \quad x_0^* \in X_0^*.$$
\end{lemma}

\textbf{Proof.}
We consider the linear mapping $\Lambda y(\cdot) = y(\cdot)-a$,
which maps the space $Y$ onto $X_0$.
Since any $y(\cdot) \in Y$ possesses the representation $y(t)=x_0(t)+a$ with $x_0(\cdot) \in X_0$, $a \in \R^n$,
the inequality $\|y(\cdot)-a\|_X=\|x_0(\cdot)\|_{X_0}\leq \|x_0(\cdot)\|_{X_0}+\|a\|= \|y(\cdot)\|_X$ holds on $Y$.
Therefore, $\Lambda$ is continuous.
Moreover, $Ker\,\Lambda=\{ 0\}$.
Let $y^* \in Y^*$.
By Theorem \ref{LemmaClosedRange},
there exists a functional $x_0^* \in X_0^*$ with $y^* = \Lambda^* x_0^*$.
This shows
$\langle y^*,y(\cdot) \rangle = \langle x_0^*, \Lambda y(\cdot) \rangle  = \langle x_0^*, y(\cdot)-a \rangle.$
The uniqueness of this representation can be verified directly. \hfill $\square$
%
\section{Spaces of Continuous Functions Converging at Infinity} \label{SectionContinuous}
We consider the space $C_{\lim}(\R_+,\R^n)$ of continuous vector-functions converging at infinity:
\begin{eqnarray*}
C_{\lim}(\R_+,\R^n) &=& \{x(\cdot) \in C(\R_+,\R^n) \,|\, \lim_{t \to \infty} x(t)= a \in \R^n \mbox{ exists}\} \\
                    &=& \{x(\cdot): \R_+ \to \R^n \,|\, x(t)=x_0(t)+a,\;x_0(\cdot) \in C_0(\R_+,\R^n),\;a \in \R^n\}.
\end{eqnarray*}
In the considerations in the previous section we equipped the space $X$ with a $p-$norm.
But in the framework of the space $C_{\lim}(\R_+,\R^n)$ the norms $\|x(\cdot)\|_\infty$
and $\|x(\cdot)\|_X=\|x_0(\cdot)\|_\infty + \|a\|$ differ.
Therefore, we have to show the equivalence of both norms on $C_{\lim}(\R_+,\R^n)$.

\begin{lemma}
On $C_{\lim}(\R_+,\R^n)$ the norms $\|x(\cdot)\|_\infty$ and $\|x(\cdot)\|_X=\|x_0(\cdot)\|_\infty + \|a\|$
are equivalent.
\end{lemma}

\textbf{Proof.} By definition,
any $x(\cdot) \in C_{\lim}(\R_+,\R^n)$ possesses the unique representation $x(\cdot)=x_0(\cdot)+a$
with $x_0(\cdot) \in C_0(\R_+,\R^n)$ and $a \in \R^n$.
Since $x(t) \to a$ as $t \to \infty$, the inequality $\|x(\cdot)\|_\infty \geq \|a\|$ holds.
Let $\|x_0(\cdot)\|_\infty \leq \|a\|$.
Then we obtain immediately $\|x(\cdot)\|_X \leq 2\|x(\cdot)\|_\infty$.
In the case $\|a\| \leq \|x_0(\cdot)\|_\infty$ there exists $\lambda \in [0,1]$ with $\|a\| = \lambda \|x_0(\cdot)\|_\infty$.
It follows
$$\|x(\cdot)\|_\infty = \|x_0(\cdot)+a\|_\infty \geq \max \{\|x_0(\cdot)\|_\infty - \|a\|,\|a\|\}
  \geq \min_{\lambda \in [0,1]} \max\{1-\lambda,\lambda\}\cdot \|x_0(\cdot)\|_\infty = \frac{1}{2} \|x_0(\cdot)\|_\infty.$$
Consequently, the inequalities $\|x(\cdot)\|_\infty \leq \|x(\cdot)\|_X \leq 3\|x(\cdot)\|_\infty$ hold on
$C_{\lim}(\R_+,\R^n)$. \hfill $\square$ \\

In the following we make use of the notation $x(\infty)=\lim\limits_{t \to \infty} x(t)$
instead of the vector $a \in \R^n$.

\begin{lemma}[Riesz's Representation Theorem on $C_{\lim}(\R_+,\R^n)$] \label{CorollaryDarstellung1}
Any continuous linear functional $x^*$ on $C_{\lim}(\R_+,\R^n)$ can be uniquely represented in the form
$$\langle x^*(\cdot),x(\cdot) \rangle = \int_0^\infty \langle x(t)-x(\infty), d\mu(t) \rangle
              +\alpha^T x(\infty) = \int_0^\infty \langle x_0(t), d\mu(t) \rangle +\alpha^T x(\infty),$$
where $\mu=(\mu_1,...,\mu_n)$ is a vector of 
measures $\mu_1,...,\mu_n \in \mathscr{M}(\R_+)$
(cf. \cite{Rudin}) and $\alpha \in \R^n$.
\end{lemma}

\textbf{Proof.}
The representation formula follows directly from Lemma \ref{LemmaRepresentation1}
with $X_0=C_0(\R_+,\R^n)$. \hfill $\square$ \\

In Lemma \ref{CorollaryDarstellung1} let
$\nu=\mu$, $\displaystyle \nu_\infty =\alpha^T-\int_0^\infty d\mu(t)$ and $\tilde{\mu}=\nu+\nu_\infty$.
Then the Riesz's representation theorem on $C_{\lim}(\R_+,\R^n)$ becomes the equivalent form:

\begin{lemma}
Any continuous linear functional $x^*$ on $C_{\lim}(\R_+,\R^n)$ can be uniquely represented by
$$\langle x^*(\cdot),x(\cdot) \rangle = \int_0^\infty \langle x(t), d\mu(t) \rangle,$$
where $\mu=(\mu_1,...,\mu_n)$ is a vector of signed regular Borel measures $\mu_1,...,\mu_n \in \mathscr{M}(\overline{\R}_+)$.
\end{lemma}

Now we consider the space $C_{\lim}(\R,\R^n)$, which we define by
$$C_{\lim}(\R,\R^n) = \{x(\cdot) \in C(\R,\R^n) \,|\, \lim_{t \to \infty} x(t)= x(\infty), \;
 \lim_{t \to -\infty} x(t)= x(-\infty) \mbox{ exist}\}.$$
Since we have to take into account both limits $x(-\infty)$ and $x(\infty)$,
the previous approach to the space $C_{\lim}(\R_+,\R^n)$ is not applicable.
Therefore, we will make use of the splitting argument.

\begin{lemma}[Riesz's Representation Theorem on $C_{\lim}(\R,\R^n)$] \label{CorollaryDarstellung3}
Any continuous linear functional $x^*$ on $C_{\lim}(\R,\R^n)$ can be uniquely represented in the form
$$\langle x^*(\cdot),x(\cdot) \rangle =
  \int_{-\infty}^\infty \langle x(t), d\mu(t) \rangle+\alpha_1^T x(-\infty)+\alpha_2^T x(\infty),$$
where $\mu$ is a vector of signed regular Borel measures on $\R$ and $\alpha_1,\alpha_2 \in \R^n$,
or equivalently in the form
$$\langle x^*(\cdot),x(\cdot) \rangle = \int_{-\infty}^\infty \langle x(t), d\tilde{\mu}(t) \rangle,$$
where $\tilde{\mu}=(\tilde{\mu}_1,...,\tilde{\mu}_n)$ is a vector of signed regular Borel measures
$\tilde{\mu}_1,...,\tilde{\mu}_n \in \mathscr{M}(\overline{\R})$.
\end{lemma}

\textbf{Proof.}
Let $Y \subset C_{\lim}(\R_+,\R^n)$ be the subspace
$Y = \{y(\cdot) \in C_{\lim}(\R_+,\R^n) \,|\, y(0) =0 \}$.
We split $x(\cdot) \in C_{\lim}(\R,\R^n)$ into $x_1(\cdot), x_2(\cdot)$ by
$\big(x_1(t), x_2(t)\big)=\big(x(-t)-x(0), x(t)-x(0)\big)$ with $t \in \R_+$.
The operator $\Lambda x(\cdot)=\big(x_1(\cdot), x_2(\cdot)\big)$ maps the space $C_{\lim}(\R,\R^n)$ onto
the product $Y \times Y$.
The mapping $\Lambda$ is linear and continuous,
and $Ker\,\Lambda$ consists of all constant functions $x(\cdot) \in C_{\lim}(\R,\R^n)$. \\
As in the proof of Lemma \ref{CorollaryDarstellung1} we define the vector $\alpha$ and
consider $x_1^* \in C^*_{\lim}(\R,\R^n)$ defined by the formula
$\langle x_1^*,x(\cdot) \rangle = \langle x^*,x(\cdot) \rangle - \alpha^T x(\infty)$.
Obviously, $x_1^* \in (Ker\,\Lambda)^\perp$.
By Theorem \ref{LemmaClosedRange},
there exists a functional $y^*=(y_1^*,y_2^*) \in Y^* \times Y^*$ with $x_1^* = \Lambda^* y^*$.
That means the equation
$$\langle x_1^*,x(\cdot) \rangle = \langle \Lambda^* y^*,x(\cdot) \rangle = \langle y^*,\Lambda x(\cdot)\rangle
  =\langle y_1^*,x_1(\cdot)\rangle + \langle y_2^*,x_2(\cdot)\rangle$$
holds for all $x(\cdot) \in C_{\lim}(\R,\R^n)$.
The functions $x_1(\cdot), x_2(\cdot) \in Y$ possessing the limits
$$\lim_{t \to \infty} x_1(t) = \lim_{t \to \infty} x(-t)-x(0) =x(-\infty)-x(0), \quad
  \lim_{t \to \infty} x_2(t) = \lim_{t \to \infty} x(t)-x(0) =x(\infty)-x(0).$$
Then Lemma \ref{LemmaRepresentation2} yields the representations
\begin{eqnarray*}
\langle y_1^*,x_1(\cdot)\rangle
&=& \int_0^\infty \langle x_1(t)-x_1(\infty), d\mu_1(t) \rangle
    = \int_0^\infty \langle x(-t)-x(-\infty), d\mu_1(t) \rangle \\
\langle y_2^*,x_2(\cdot)\rangle
&=& \int_0^\infty \langle x_2(t)-x_2(\infty), d\mu_2(t) \rangle
    = \int_0^\infty \langle x(t)-x(\infty), d\mu_2(t) \rangle,
\end{eqnarray*}
where $\mu_1,\mu_2$ are vectors of signed regular Borel measures on $\R_+$.
By definition of $x_1^*$ we obtain
$$\langle x^*,x(\cdot)\rangle =
  \int_0^\infty \langle x(-t)-x(-\infty), d\mu_1(t) \rangle
  + \int_0^\infty \langle x(t)-x(\infty), d\mu_2(t) \rangle + \alpha^T x(\infty).$$
Finally, we determine $\alpha_1$, $\alpha_2$, $\tilde{\mu}_1$ and $\mu$ by
$$\alpha_1= -\int_0^\infty d\mu_1(t), \quad
  \alpha_2= \alpha-\int_0^{\infty} d\mu_2(t), \quad \tilde{\mu}_1(t)=\mu_1(-t), \quad \mu=\tilde{\mu}_1+\mu_2.$$
Then the representation formula in Lemma \ref{CorollaryDarstellung3} is shown. \hfill $\square$
\section{Lebesgue Spaces of Functions Converging at Infinity} \label{SectionLebesgue}
Throughout this section let $1 \leq p < \infty$
and let $q$ determined by $1/p+1/q=1$ if $p>1$ and $q=\infty$ if $p=1$.
Since the representation of $x \in X$ in Section \ref{SectionConcepts} has to be unique,
the parameter $p$ must be restricted to $p < \infty$.
For an introduction of the Lebesgue spaces $L_p(\R_+,\R^n)$ cf. \cite{Werner}. \\

We consider the Lebesgue space $L_{p,\lim}(\R_+,\R^n)$ defined by
$$L_{p,\lim}(\R_+,\R^n)= \{x(\cdot): \R_+ \to \R^n
           \,|\, x(t)=x_0(t)+a,\; x_0(\cdot) \in L_p(\R_+,\R^n),\; a \in \R^n\}.$$
The space $L_{p,\lim}(\R_+,\R^n)$ is equipped with the norm $\|x(\cdot)\|_{L_{p,\lim}}^p=\|x_0(\cdot)\|_{L_p}^p+\|a\|^p$.
Furthermore, any $x(\cdot) \in L_{p,\lim}(\R_+,\R^n)$ possesses at infinity a $\lambda-$limit.
We write $x(\infty)=\llim\limits_{t \to \infty} x(t)=a$.

\begin{lemma} \label{LemmaLebesgue}
Let $1 \leq p < \infty$. 
Any $x^* \in L^*_{p,\lim}(\R_+,\R^n)$ can be uniquely represented in the form
$$\langle x^*(\cdot),x(\cdot) \rangle
   = \int_0^\infty \langle y(t), x(t)-x(\infty) \rangle \, dt + \alpha^T x(\infty), \qquad
  y(\cdot) \in L_q(\R_+,\R^n),\; \alpha \in \R^n.$$
\end{lemma}

\begin{pf}
The representation formula follows directly from Lemma \ref{LemmaRepresentation1} with $X_0=L_p(\R_+,\R^n)$. \hfill $\square$
\end{pf}

Now, we consider the Lebesgue space $L_{p,\lim}(\R,\R^n)$, which we define by
\begin{eqnarray*}
&& L_{p,\lim}(\R,\R^n)= \{x(\cdot): \R \to \R^n
           \,|\, x(t)=x_0(t)+a_1 \mbox{ on } \R_-,\;
                 x(t)=x_0(t)+a_2 \mbox{ on } \R_+, \\
&& \hspace*{5cm} x_0(\cdot)\in L_p(\R,\R^n),\; a_1,a_2 \in \R^n\}.
\end{eqnarray*}
$L_{p,\lim}(\R,\R^n)$ is equipped with the norm $\|x(\cdot)\|_{L_{p,\lim}}^p=\|x(\cdot)\|_{L_p}^p+\|a_1\|^p+\|a_2\|^p$.
Moreover, we determine the $\lambda-$limits $x(-\infty)=a_1$ and $x(\infty)=a_2$.

\begin{lemma} \label{LemmaLebesgue2}
Let $1 \leq p < \infty$. 
Any $x^* \in L^*_{p,\lim}(\R,\R^n)$ can be uniquely represented in the form
$$\langle x^*(\cdot),x(\cdot) \rangle =
  \int_{-\infty}^\infty \langle y(t), x(t)\rangle \, dt
  + \alpha_1^T x(-\infty)+ \alpha_2^T x(\infty), \qquad y(\cdot) \in L_q(\R,\R^n),\; \alpha_1,\alpha_2 \in \R^n.$$
\end{lemma}

\textbf{Proof.}
Any $x(\cdot) \in L_{p,\lim}(\R,\R^n)$ can be split into $x_1(\cdot), x_2(\cdot) \in L_{p,\lim}(\R_+,\R^n)$ by
$$\Lambda x(\cdot)=\big(x_1(\cdot), x_2(\cdot)\big), \qquad \big(x_1(t),x_2(t)\big)=\big(x(-t), x(t)\big), \quad t \in \R_+.$$
The operator $\Lambda$ maps the space $L_{p,\lim}(\R,\R^n)$ onto $L_{p,\lim}(\R_+,\R^n) \times L_{p,\lim}(\R_+,\R^n)$.
Moreover, $\Lambda$ is linear and continuous, and $Ker\,\Lambda =\{0\}$.
By Theorem \ref{LemmaClosedRange},
any $x^* \in L^*_{p,\lim}(\R,\R^n)$ possesses the representation
$$\langle x^*,x(\cdot) \rangle
  = \langle \Lambda^* y^*,x(\cdot) \rangle = \langle y^*,\Lambda x(\cdot)\rangle
  =\langle y_1^*,x_1(\cdot)\rangle + \langle y_2^*,x_2(\cdot)\rangle, \quad y^*_1,y_2^* \in L^*_{p,\lim}(\R_+,\R^n).$$
That means, there are unique $\alpha_1,\alpha_2 \in \R^n$ and $y_1(\cdot),y_2(\cdot) \in L_q(\R_+,\R^n)$ with
\begin{eqnarray*}
&&\hspace*{-8mm}\langle y_1^*,x_1(\cdot) \rangle
    = \int_0^\infty \langle y_1(t), x_1(t)-x_1(\infty) \rangle \, dt + \alpha_1^T x_1(\infty)
    =\int_0^\infty \langle y_1(t), x(-t)-x(-\infty) \rangle \, dt + \alpha_1^T x(-\infty), \\
&&\hspace*{-8mm}\langle y_2^*,x_2(\cdot) \rangle
    = \int_0^\infty \langle y_2(t), x_2(t)-x_2(\infty) \rangle \, dt + \alpha_2^T x_2(\infty)
    =\int_0^\infty \langle y_2(t), x(t)-x(\infty) \rangle \, dt + \alpha_2^T x(\infty).
\end{eqnarray*}
Finally, the representation formula follows
by the determination of $y(\cdot)$ by $y(t)=y_1(-t)$ on $\R_-$ and $y(t)=y_2(t)$ on $\R_+$,
and of $\tilde{\alpha}_1,\tilde{\alpha}_2$ by
$\tilde{\alpha}_1=\displaystyle \alpha_1-\int_{-\infty}^0 y(t) \, dt$,
$\tilde{\alpha}_2=\displaystyle \alpha_2-\int_0^{\infty} y(t)\,dt$. \hfill $\square$
\section{Remarks} \label{SectionRemarks}
\textbf{Spaces of continuous functions:} 
The introduction of the space $C_{\lim}(\R_+,\R^n)$ were motivated by degenerations in the convex optimization on $\R_+$.
In this context we note:
\begin{enumerate}
\item The interior of the cone $\mathscr{K}=\{ z(\cdot) \in C_{\lim}(\R_+,\R) \,|\, z(t)\leq 0 \mbox{ for all } t \in \R_+\}$
      is non-empty and containing any $x(\cdot) \in C_{\lim}(\R_+,\R)$ with
      $\max\limits_{t \in \overline{\R}_+} x(t) <0$.      
\item We consider $f\big(x(\cdot)\big)=\max\limits_{t \in \overline{\R}_+} x(t)$.
      By definition of $\partial f(0)$ we obtain in the space $C_{\lim}(\R_+,\R)$:
      $$\max_{t \in \overline{\R}_+} x(t) \geq \int_0^\infty x(t) \, d\mu(t)
        \geq - \max_{t \in \overline{\R}_+} \big( -x(t)\big)
        = \min_{t \in \overline{\R}_+} x(t).$$ 
      This shows, that the subdifferential of $f$ consists in $x(\cdot)=0$ of those and only those non-negative Borel measures
      $\mu \in \mathscr{M}(\overline{\R}_+)$, which satisfy $\|\mu\| = 1$.
\end{enumerate}

\textbf{Lebesgue spaces:} The introduction of the spaces $L_{p,\lim}(\R_+,\R^n)$ and the Lemma \ref{LemmaLebesgue} deliver:
\begin{enumerate}
\item[(a)] If $1<p<\infty$, the spaces $L_{p,\lim}(\R_+,\R^n)$ are reflexive Banach spaces
           and $y(\cdot) \in L_q(\R_+,\R^n)$ and $\alpha \in \R^n$ in Lemma \ref{LemmaLebesgue}
           determine $z(\cdot)=y(\cdot)+\alpha \in L_{q,\lim}(\R_+,\R^n)$.
\item[(b)] The space $L_{2,\lim}(\R_+,\R^n)$ is a separable Hilbert space.
\item[(c)] Any $x(\cdot) \in L_{2,\lim}(\R_+,\R^n)$ can be uniquely represented in the form $x(\cdot)=x_0(\cdot)+a$ with
          $\big(x_0(\cdot),a\big) \in L_2(\R_+,\R^n) \times \R^n$.
          Let $x(\cdot),y(\cdot) \in L_{2,\lim}(\R_+,\R^n)$ with $x(\cdot)=x_0(\cdot)+a$, $y(\cdot)=y_0(\cdot)+b$.
          Then the inner product in the space $L_{2,\lim}(\R_+,\R^n)$ has the form
          $$\langle x(\cdot),y(\cdot)\rangle_{L_{2,\lim}}=
            \big\langle \big(x_0(\cdot),a\big),\big(y_0(\cdot),b\big)\big\rangle_{L_{2,\lim}}=
            \langle x_0(\cdot),y_0(\cdot)\rangle_{L_2}+a^T b.$$
\item[(d)] Let $\{\varphi_0(\cdot),\varphi_1(\cdot),...\}$ be an orthonormal basis in $L_2(\R_+,\R)$.
           The introduction of $L_{2,\lim}(\R_+,\R)$ leads to the system
          $\{(0,1),\,\big(\varphi_0(\cdot),0\big),\,\big(\varphi_1(\cdot),0\big),\,...\}$,
         which is an orthonormal basis in $L_{2,\lim}(\R_+,\R)$.
\item[(e)] In the application of Hilbert space methods in the frame work of the space $L_{2,\lim}(\R,\R)$
           the splitting concept suggests
           to split $x(\cdot) \in L_{2,\lim}(\R,\R)$ into $\big(x_1(\cdot),x_2(\cdot)\big)$ and to consider
            the Fourier series expansions of $x_1(\cdot),x_2(\cdot) \in L_{2,\lim}(\R_+,\R)$ separately.
\end{enumerate}

\textbf{Sequence spaces:} Let $1\leq p < \infty$. 
We determine the sequence space $\ell_{p,\lim}$ by
$$X=\ell_{p,\lim}= \{x=(x_n)_{n \in \N} \,|\, x_n=x^0_n+a, x^0=(x^0_n)_{n \in \N} \in \ell_p, a \in \R\},
  \quad \|x\|^p_{\ell_{p,\lim}}=\|x^0\|^p_{\ell_p}+|a|^p.$$
The elements of the spaces $\ell_{p,\lim}$ are convergent sequences.

\begin{lemma} \label{LemmaSequence}
Let $1 \leq p < \infty$. 
Any $x^* \in \ell_{p,\lim}^*$ can be uniquely represented in form
$$\langle x^*,x \rangle = \sum_{n \in N} y_n (x_n-a)+ \alpha a = \sum_{n \in N} y_n x^0_n+ \alpha a,
\qquad y \in \ell_q,\; \alpha \in \R.$$
\end{lemma}

\textbf{Proof.}
The first scheme in Section \ref{SectionConcepts} holds if we consider mappings from $\N$ to $\R^n$.
Therefore,
the representation formula follows directly from Lemma \ref{LemmaRepresentation1} with $X_0=\ell_p$. \hfill $\square$ \\

\textbf{Sobolev spaces:} Let $1 \leq p < \infty$
and let $q$ determined by $1/p+1/q=1$ if $p>1$ and $q=\infty$ if $p=1$.
For an introduction of Sobolev spaces $W^1_p(\R_+,\R^n)$ we refer to \cite{Werner}. \\

We consider the Sobolev space $W^1_{p,\lim}(\R_+,\R^n)$ defined by
$$W^1_{p,\lim}(\R_+,\R^n)= \{x(\cdot): \R_+ \to \R^n
           \,|\, x(t)=x_0(t)+a,\; x_0(\cdot) \in W^1_p(\R_+,\R^n),\; a \in \R^n\}.$$
We equip the space $W^1_{p,\lim}(\R_+,\R^n)$ with the norm
$\|x(\cdot)\|_{W^1_{p,\lim}}^p=\|x(\cdot)\|_{W^1_p}^p+\|a\|^p$.
Furthermore, $x(\cdot) \in W^1_{p,\lim}(\R_+,\R^n)$ possesses at infinity the $\lambda-$limit
$x(\infty)=\llim\limits_{t \to \infty} x(t)=a$.

\begin{lemma} \label{LemmaSobolev2}
Let $1 \leq p < \infty$. 
Any $x^* \in \big(W^1_{p,\lim}(\R_+,\R^n)\big)^*$ can be uniquely represented in the form
$$\langle x^*(\cdot),x(\cdot) \rangle
   = \langle x_0^*(\cdot),x(\cdot)-x(\infty) \rangle + \alpha^T x(\infty),
   \qquad x_0^* \in \big(W^1_p(\R_+,\R^n)\big)^*, \; \alpha \in \R^n.$$
\end{lemma}

\textbf{Proof.}
The assertion follows directly by Lemma \ref{LemmaRepresentation1}. \hfill $\square$

%

\end{document}